\documentclass[10pt]{article}
\usepackage{graphicx}
\usepackage{amssymb,amsmath}
\usepackage{float}
\usepackage{tcolorbox}
\usepackage{amsthm}
\usepackage{graphicx}
\usepackage{hyperref}
\usepackage{multirow}
\usepackage{tabularx}
\usepackage{rotating}

\usepackage[english]{babel}
\usepackage{array}
\usepackage{multirow}
\usepackage{siunitx}

\hypersetup{
	colorlinks=true,
	linkcolor=blue,
	filecolor=magenta,      
	urlcolor=cyan}

\begin{document}

\title{Echoes of the hexagon: remnants of hexagonal packing inside regular polygons}

\author{
  Paolo Amore \\
  \small Facultad de Ciencias, CUICBAS, Universidad de Colima,\\
  \small Bernal D\'{i}az del Castillo 340, Colima, Colima, Mexico \\
  \texttt{paolo@ucol.mx}
  \and
  Mauricio Carrizalez \\
  \small Facultad de Ciencias, Universidad de Colima,\\
  \small Bernal D\'{i}az del Castillo 340, Colima, Colima, Mexico \\
  \texttt{mcarrizalez@ucol.mx }
  \and
  Ulises Zarate \\
  \small Facultad de Ciencias, Universidad de Colima,\\
  \small Bernal D\'{i}az del Castillo 340, Colima, Colima, Mexico \\
  \texttt{uzarate@ucol.mx}
}

\maketitle

\begin{abstract}
Based on numerical simulations that we have carried out, we provide evidence that for regular polygons with  $\sigma= 6j$ sides (with $j=2,3,\dots$), $N(k)=3 k (k+1)+1$ (with $k=1,2,\dots$)  congruent disks  of appropriate size can be nicely packed inside these polygons in highly symmetrical configurations which apparently have maximal density for $N$ sufficiently small.  These configurations are invariant under rotations of $\pi/3$ and are closely related to the configurations with perfect hexagonal packing in the regular hexagon and to the configurations with {\sl curved hexagonal packing}  (CHP) in the circle found long time ago by Graham and Lubachevsky~\cite{Graham97}.

At the basis of our explorations are the algorithms that we have devised, which are very efficient in producing the CHP and more general configurations inside regular polygons. We have used these algorithms to generate a large number of CHP configurations for different regular polygons and  numbers of disks; a careful study of these results has made possible to fully characterize the general properties of the CHP configurations and to devise a {\sl deterministic} algorithm that completely ensembles a given CHP configuration once an appropriate input ("DNA") is specified. Our analysis shows that the number of CHP configurations for a given $N$ is highly degenerate in the packing fraction and it can be explicitly calculated in terms of $k$ (number of shells), of the building block of the DNA itself and of the number of vertices in the fundamental domain (because of the symmetry we work in $1/6$ of the whole domain). With the help of our deterministic algorithm we are able to build {\sl all} the CHP configurations for a polygon 
with $k$ shells.

We have also found few examples of these symmetric configurations with deformed hexagonal packing in certain polygons with $3 (2j+1)$ sides ($j=1,2,\dots$) as well as  configurations which have the correct symmetry and contain a perfect hexagonal packing in a smaller portion of the domain. 

\end{abstract}

\maketitle

\section{Introduction}
\label{sec:intro}

According to Thue's theorem the maximal density (packing fraction) that can be achieved in the packing of congruent disks in the plane 
is $\rho = \pi/\sqrt{12}$ and it corresponds to forming a regular structure in which each disk is at the center of a regular hexagon with other $6$ disks located at the vertices of the hexagon. 

However, if one considers a finite region instead of the whole plane,  hexagonal packing is generally not possible because of the geometric frustration introduced  by the border. An exception to this rule occurs when the container is the regular hexagon, with $N$ disks and $N = h(k) \equiv 3 k (k+1)+1$ 
and $k \geq 1$ (see for example ref.~\cite{Amore22}), or the equilateral triangle, with $N = t(k) \equiv k (k+1)/2$ (see for example refs.~\cite{Graham04} and \cite{Amore22}).
 
Long time ago Graham and Lubachevsky \cite{Graham97} while studying the packing of congruent disks inside a circle found that the configurations with 
largest density for $N = h(k)$ with $k\leq 5$, correspond to a curved hexagonal packing (CHP) pattern which is invariant under rotation of $60^0$ (see for example Fig. 1.1 of \cite{Graham97}). The conclusions reached by those authors are based on large scale numerical computations, where the densest configurations obtained
for the cases above were always CHP. Starting at $N= 127$ ($k=6$) these configurations are no longer optimal, since denser configurations with no symmetry were obtained. 

While studying the packing of congruent disks inside regular polygons of different number of sides $\sigma$ in ref.~\cite{Amore22} one of us realized that the phenomenon first observed by  Graham and Lubachevsky for the circle also occurs inside regular polygons, particularly if $\sigma$ is a integer multiple of $6$ (the surprise can be mitigated if one realizes that the circle is just a regular polygon with $\sigma=\infty$). 
The purpose of the present paper is to explore this uncharted territory and provide a full characterization of the properties of these peculiar (and beautiful) configurations.


The paper is organized as follows: in section \ref{sec:algo} we briefly describe the numerical algorithms (the reader may refer to \cite{Amore22} for a more detailed  discussion); in section \ref{sec:chp} we establish the general properties of the CHP configurations and in section \ref{sec_detalg} set up a {\sl deterministic} algorithm that is capable to fully build {\sl all} the configurations for a given polygon and a given number of shells (i.e. $\sigma$ and $k$); in section \ref{sec:res} we discuss the results obtained for different cases; in section \ref{sec:break} we explore the use of CHP configurations with 
a large number of disks (which are certainly not optimal) to achieve much denser disordered configurations, using algorithm 2 of \cite{Amore22}; 
 finally in section \ref{sec:concl} we draw our conclusions and briefly  suggest future directions of work.

\section{Numerical algorithms}
\label{sec:algo}

The algorithms that we use in this paper have been recently developed in ref.~\cite{Amore22}; in particular the first algorithm of ref.~\cite{Amore22} is the evolution of an algorithm originally proposed by Nurmela and Östergård (N\"O) long time ago for the square~\cite{Nurmela97}. Recently this algorithm has been  modified and improved in ref.~\cite{Amore22b}, by introducing border repulsion \footnote{Notice also that the algorithms of \cite{Nurmela97,Amore22b} are limited to square container, whereas the algorithms of \cite{Amore22} can deal with any regular polygon with an arbitrary number of sides.}.

We will limit ourselves to sketch briefly the basic idea of the algorithm here (the reader will find all the needed details in \cite{Amore22}) and then describe how these algorithms (particularly algorithm 2) can be used to generate curved hexagonal packing inside  regular polygons with $6j$ sides.

Following ref.~\cite{Amore22} we parametrize the coordinates of a point inside a regular polygon with $\sigma$ sides  as
\begin{equation} 
\begin{split} 
x(t,u,\sigma) &= \sin^2 t \ X(u,0,\sigma) \\
y(t,u,\sigma) &= \sin^2 t \ Y(u,0,\sigma)   \ , 
\label{eq_poly_par}
\end{split}
\end{equation}
where 
\begin{equation} 
\begin{split} 
X(u,\delta,\sigma) &\equiv \Gamma(u,\delta,\sigma)  \cos (u) \\
Y(u,\delta,\sigma) &\equiv \Gamma(u,\delta,\sigma)  \sin (u) \ ,
\label{eq_poly_par2}
\end{split}
\end{equation}
and
\begin{equation}
\Gamma(u,\delta,\sigma) \equiv \left(\delta +\cos \left(\frac{\pi }{\sigma}\right)\right) \sec
\left(\frac{\pi }{\sigma}-(u \bmod \frac{2 \pi }{\sigma})\right) \ .
\end{equation}

In particular $P = (X(u,0,\sigma),Y(u,0,\sigma))$ is a point on the border of the polygon.

Each point is considered to be the center of a disk of radius $R$: the disk are allowed to be in contact with other disks and with the container, but not to overlap.
In this way 
\begin{equation}
R = \min_{1 \leq i \neq j \leq N } \left\{ r_{ij}\right\}
\end{equation}
and one wants to dispose the points in the domain such that $R$ is maximal (in this formulation packing becomes a {\sl minimax} problem). 
It is important to observe that there are two polygons: an inner polygon, represented by the points $P$ above, which corresponds to the region of the plane where the centers of the disks are  allowed to move, and an outer  polygon which corresponds to the actual container (see for example Fig.1 of \cite{Amore22}).

In the formulation of N\"O the points are regarded as charges which interact with a potential $(\lambda/r^2)^s$ and are initially randomly distributed inside the polygon (in their case the polygon is only a square). The value of $s$ determines the range of the potential: a small $s$ corresponds to a long-range potential, where a single point interacts with many other points, even at large distance; a large $s$ corresponds to a short range potential which  is very large for $0 \leq r \leq \sqrt{\lambda}$ and dies off quite rapidly for $r > \sqrt{\lambda}$. The parameter $\lambda$ is clearly related to the radius of the disks and, for $s \rightarrow \infty$, one obtains the expected contact interaction between rigid disks.

In the algorithm 1 one starts with a random initial configuration and minimizes the total energy of the system for a given value of $s$ (at this stage $s$ should be not too large, so that the long range interaction can easily establish an order over all the domain). Once the new arrangement is found, the value of $s$ is increased and the minimization is repeated. This process is repeated many times until $s$ reaches very large values ($s \approx 10^8 - 10^9$ in our case). The final configuration obtained in this way corresponds to a circle packing in the polygon, but it is not necessarily optimal (in fact the algorithm does not directly maximizes the packing density). For this reason, the process has to be repeated a large number of times and the best configuration obtained can then be regarded as a candidate for global maximum of the packing fraction~\footnote{Using their algorithm Nurmela and Östergård were able to prove in \cite{Nurmela97} that the optimal configuration of $49$ congruent disks inside the square does not correspond to a square packing.}. 

Algorithm 2 was introduced originally in ref.~\cite{Amore22b} for the square and later adapted to general regular polygons in \cite{Amore22}.
In this case the initial configuration is not random, but it could be a packing configuration obtained using algorithm 1 (or any other packing algorithm). The basic idea of the algorithm come from Physics: experiments show that the packing density of objects inside a container is larger if the container in which the objects are poured is shaken~\cite{Pouliquen97,Kudrolli2010}. The shaking in our case is introduced by random perturbations of the positions of the disks, followed by a minimization of the total interaction energy, for increasing values of $s$. If the density at the end of the process has increased, the new configuration is retained, otherwise one keeps the original configuration. The process is repeated a large number of times.  This algorithm can be very successful at improving even very dense configurations of large numbers of disks (such as the case of $3000$ disks in the square considered in \cite{Amore22b} ).

Finally, a third algorithm (refinement algorithm) can also be used to "refine" the positions of the disks to high precision 
(see the discussion in Ref.~\cite{Amore22}).

In this paper we use the three algorithms described above, particularly algorithm 2, but with some important adjustments. Our goal, in fact, is to obtain configurations of curved hexagonal packing (CHP) inside a polygon. These configurations are related to the configurations of hexagonal packing inside the hexagon, which occur for special values $N(k) = 3 k (k+1)+1$ and $k=1,2,\dots$. 

\begin{figure}
\begin{center}
\bigskip\bigskip\bigskip
\includegraphics[width=6cm]{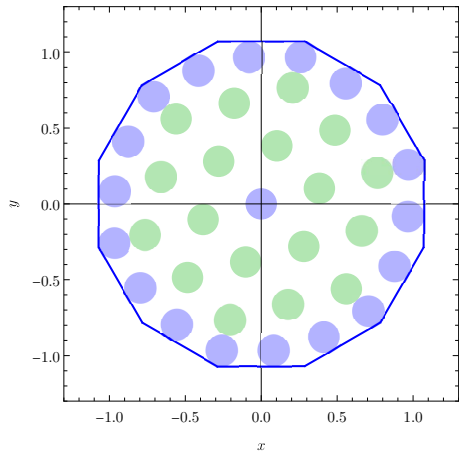}\hspace{1cm}
\bigskip
\caption{Initial configuration of $37$ disks inside a dodecagon: the blue disks are already at the exact positions, whereas the green disks correspond to hexagonal packing rotated by an arbitrary angle and rescaled to fit the domain.}
\label{Fig_1}
\end{center}
\end{figure}

CHP configurations, therefore, must be invariant under rotations of multiples of $\pi/3$: for this reason it is natural to look for CHP in 
polygonal domains with $\sigma = 6j$ sides, with $j=2,3,\dots$. 

A second property of hexagonal packings which should also be retained by CHP is  their "shell" structure: in a configuration with  $N(k)$ disks there are $6k$ outer  disks that are at contact with the container. The density of a CHP can therefore be calculated by determining the positions of $k$ points uniformly distributed on the border for $0 \leq u < \pi/3$. The remaining points on the border can be obtained by using the symmetry under  rotations of multiples of $\pi/3$. Another special point is the origin, which always corresponds to the position of the central disk.

The algorithm may be much more effective if the initial configuration is chosen properly: we consider a configuration where the $6k$ border disks and the center disk are already located at their exact positions and the remaining $3 k (k-1)$ points correspond to hexagonal packing inside the domain, rotated of some arbitrary angle and suitably rescaled (see Fig.~\ref{Fig_1}).

It is convenient to take advantage of the fact that there are $6k+1$ points already at their exact positions; clearly these points should not be moved by the algorithm! This condition can be enforced by setting to zero the components of the gradient corresponding to the exact points. The algorithm modified in this way will be faster and also much more effective in finding CHP configurations, since the internal points are now guided to their final position in a much more constrained way. As a matter of fact the probability of generating a CHP in this way with a single trial is considerably larger that using algorithm 1 a large number of times, or even algorithm 2 with a generic initial configuration. 

Even so, finding all or most of the degenerate CHP configurations for a given $N(k)$ is a challenging and time--consuming process, depending on how large $N$ is; it is possible to drastically improve the performance of our algorithm above in a simple, yet very effective way, with a limited computational cost.

\begin{figure}
\begin{center}
\bigskip\bigskip\bigskip
\includegraphics[width=6cm]{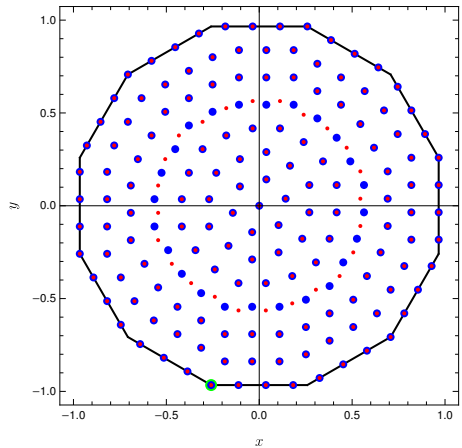}
\hspace{1cm} 
\caption{Comparison of two CHP configurations of 169 disks inside the dodecagon (the blue and red points correspond to the centers of the disks for each configuration while the green point marks the position of the fundamental vertex).  The solid line represents the dodecagon over which the centers of the border disks are distributed.}
\label{Fig_ALG}
\end{center}
\end{figure}

To illustrate how this works we consider Fig.~\ref{Fig_ALG}, where two independent CHP configurations are plotted for the case of 169 disks inside the dodecagon. The blue and red points correspond to the two cases. We see that these two configurations share a large number of points, but 
a particular shell (the fourth from the center) is rotated. As a result, if one rotates the blue points of the fourth shell by an appropriate amount
and then applies algorithm 2, the probability of generating the second configuration is very large. With this idea in mind, we have operated in the following way:
\begin{itemize}
	\item Assume that one has at least one CHP
	\item Perform several trials, where different shells (randomly chosen) are rotated by a finite (random) amount;
	\item Apply algorithm 2 to each of these cases and for each case verify whether the corresponding configuration is a CHP or not; if it is a CHP compare it with the original CHP and accept it if it is different from it;
	\item After finishing these trials, repeat the whole process, this time using one of the independent CHPs that you have found with the previous process;
\end{itemize}

We have tested this procedure over the case of $169$ disks inside the dodecagon (which does not correspond to a global maximum of the packing density)
finding a large number of configurations that we hadn't found earlier. 

\section{Curved hexagonal packing}
\label{sec:chp}

We may use the parametrization of the regular polygon in eq.~(\ref{eq_poly_par}) to obtain under which conditions
CHP is possible and what is the corresponding density of the packing. 

\begin{figure}[H]
\begin{center}
\bigskip\bigskip\bigskip
\includegraphics[width=7cm]{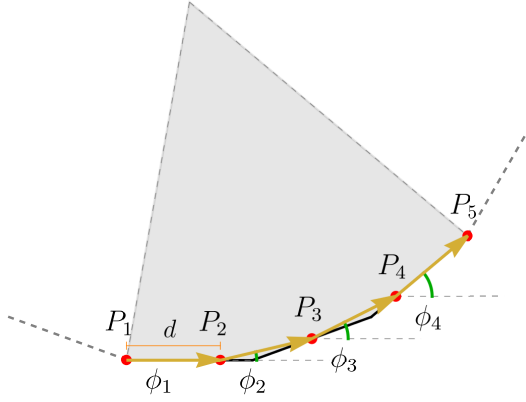}
\hspace{1cm} 
\caption{Determination of the border points  for the case of the octadecagon with $k=4$. The shaded area corresponds to the fundamental region.}
\label{Fig_P1}
\end{center}
\end{figure}

In what follows we adopt the convention that the regular polygon is oriented as shown in Fig.~\ref{Fig_ALG} and consider as fundamental vertex the point $P_1 \equiv (-\sin \frac{\pi}{\sigma},-\cos \frac{\pi}{\sigma})$ which is painted in green in the figure. Because of the symmetry of the configurations under rotations of $\pi/3$ we may restrict our considerations to the fundamental region $\frac{3\pi}{2} -\frac{\pi}{\sigma} \leq \phi \leq \frac{11\pi}{6} -\frac{\pi}{\sigma}$. For a configuration with $k$ shells, i.e. $N(k) = 3 k (k+1)+1$ disks, there  are $6k$ disks arranged on the border and therefore $k$ disks in the fundamental region. 
By convention we place the first disk at the vertex $P_1$; assuming that the disks have a diameter $d$, the next disks will be placed at the points (on the regular polygon):
\begin{equation}
\begin{split}
P_2 &= P_1 + d (\cos \varphi_1, \sin \varphi_1) \\
P_3 &= P_2 + d (\cos \varphi_2, \sin \varphi_2) \\
&\dots
\end{split}
\end{equation}
with the condition
\begin{equation}
P_{k+1} = \left(\cos \left(\pi \left(\frac{1}{\sigma}+\frac{1}{6}\right)\right),-\sin\left(\pi  \left(\frac{1}{\sigma}+\frac{1}{6}\right)\right)\right) \ .
\label{eq_condition}
\end{equation}
Here $\varphi_j$ are the angles that the segment joining the $j$ and $(j+1)$ disks forms with respect to the horizontal axis. This construction is illustrated in Fig.~\ref{Fig_P1} for the case of the octadecagon with $k=4$.

Eq.~(\ref{eq_condition}) is fulfilled only if the diameter of the disks is given by the formula
\begin{equation}
d = \frac{\sin \left(\frac{\pi }{\sigma }\right)+\cos \left(\frac{\pi}{\sigma }+\frac{\pi }{6}\right)}{\sum _{j=1}^k \cos (\phi_j)}
\end{equation} 
and the angles obey the supplementary condition
\begin{equation}
\frac{\sum _{j=1}^k \sin (\phi_j)}{\sum _{j=1}^k \cos (\phi_j)}  = \frac{4}{\tan \left(\frac{\pi }{\sigma }\right)+\sqrt{3}}-\sqrt{3} \ .
\end{equation} 

The packing fraction (density) is defined as the ratio between the area covered by the disks and the total area of the polygon and
it can be expressed as
\begin{equation}
\begin{split}
\rho(k,\sigma) &=  (3 k (k+1)+1) \frac{\pi  d^2 \cot \left(\frac{\pi }{\sigma}\right)}{\sigma  \left(d+2 \cos \left(\frac{\pi }{\sigma}\right)\right)^2} \\
&= 8 \pi  (3 k (k+1)+1) \frac{ \csc \left(\frac{2 \pi }{\sigma }\right)
   \left(\sin \left(\frac{\pi }{\sigma }\right)+\cos \left(\frac{\pi
   }{\sigma }+\frac{\pi }{6}\right)\right)^2}{\sigma  \left(4 \sum
   _{j=1}^k \cos (\phi (j))+\tan \left(\frac{\pi }{\sigma
   }\right)+\sqrt{3}\right){}^2} 
\end{split}
\label{eq_rho_chp}
\end{equation}


For the special case in which all the vertices are occupied it must occur that $6k/\sigma$ is an integer and 
\begin{equation}
\phi_j = \left\{ \begin{array}{ccc}
0 & , & j=1, \dots, 6k/\sigma \\
2\pi/\sigma & , & j=6k/\sigma+1, \dots, 12 k/\sigma \\
\dots & & \\
2\pi/6 &, &  j = k-6k/\sigma +1, \dots, k \\
\end{array} \right. \ .
\end{equation}

In this case the density reduces to
\begin{equation}
\rho(k,\sigma) = \frac{\pi  (3 k (k+1)+1) \sigma  \cot\left(\frac{\pi}{\sigma}\right)}{\left(6 k \cot\left(\frac{\pi}{\sigma}\right)+\sigma\right)^2} \ ,
\label{eq_rho_chp_full}
\end{equation}
which for the hexagon reduces to 
\begin{equation}
\begin{split}
\rho(k,6) &= \frac{6 \sqrt{3} \pi  (3 k (k+1)+1)}{\left(6 \sqrt{3} k+6\right)^2} \\
\end{split} \ .
\label{eq_rho_chp_hex}
\end{equation}

Similarly the case of the circle corresponds to using $\phi_j= (2j-1) \pi/6k$ inside eq.~(\ref{eq_rho_chp}) obtaining
\begin{equation}
\begin{split}
\rho(k,\infty) &= (3 k (k+1)+1) \frac{\sin^2 \frac{\pi}{6k}}{\left(1+\sin\frac{\pi}{6k}\right)^2} \\
\end{split} \ .
\label{eq_rho_chp_circle}
\end{equation}

From our previous analysis we see that for a regular polygon with $\sigma$ sides and $k$ shells there is a {\sl unique} 
set of values $\left\{ \varphi_1, \varphi_2, \dots, \varphi_k \right\}$ that allows curved hexagonal packing. In general, however, there will be many nonequivalent packing configurations which correspond to the same set of angles $\varphi$~\footnote{Two configurations are nonequivalent if they cannot be brought to coincide by means of rotations and reflections applied to one of them.}. These configurations will have the same density and the same arrangement of disks on the border. 
Numerical evidence for the dodecagon, for instance, suggests that the number of "degenerate" configurations grows 
as $\frac{k!}{2 (\lfloor k/2\rfloor )!^2}$~\footnote{This expression turns out to be a special case of the general formula that we have derived in this paper and the we present later.}.

\begin{figure}[H]
\begin{center}
\bigskip\bigskip\bigskip
\includegraphics[width=5cm]{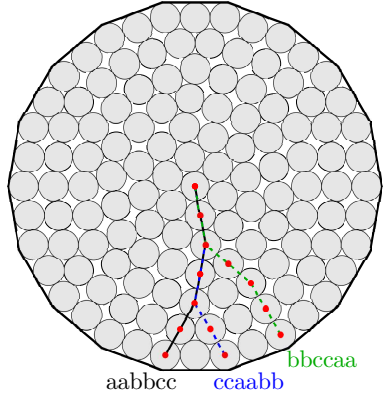}
\hspace{1cm} 
\caption{DNA of a configuration of the octadecagon with $k=6$.}
\label{Fig_DNA}
\end{center}
\end{figure}

To fully characterize a configuration we consider the vertex $P_1$ and draw the shortest path from $P_1$ passing for $k$ disks and ending in the origin (see Fig.~\ref{Fig_DNA}). The segment joining two $j^{th}$ and $(j+1)^{th}$ disks will have length $d$ and forms an angle $\xi_j$ with the horizontal axis. For each configuration, however, there will be two possible sets of angles $\xi$ associated, because of the  possibility of performing a reflection about the axis that goes through $P_1$ and the origin.

Under this reflection one has
\begin{equation}
\xi_j \rightarrow \pi-\xi_j  - \frac{2\pi}{\sigma} \hspace{1cm} j=1,2,\dots,k \ .
\label{eq_refl}
\end{equation}

Thus $\left\{ \xi_1, \dots, \xi_k \right\}$ and $\left\{ \pi-\xi_1  - \frac{2\pi}{\sigma}, \dots, \pi-\xi_k  - \frac{2\pi}{\sigma} \right\}$ represent the same configuration. This implies that $ \pi-\xi_j  - \frac{2\pi}{\sigma} = \xi_l$ for some $1 \leq l \leq k$.

To help with the classification of these configurations we adopt the convention of associating a letter to each angle, starting from the smallest angle and  proceeding in increasing order. For instance, in the dodecagon with $61$ disks ($k=2$) we only observe the values $\pi/3$ (twice) and $\pi/2$ (twice) for $\xi$ . The first angle is then represented by the letter $a$, while the second angle by the letter $b$. In this case the reflection (\ref{eq_refl}) corresponds to the interchange $a \leftrightarrow b$ and therefore $\mathbf{aabb}$ and  $\mathbf{bbaa}$ are the same configuration.  
We find convenient to use always the notation corresponding to the lowest lexicographic order (in this case $\mathbf{aabb}$).

We have carried out a large number of numerical experiments for  different $\sigma$ and $k$ and produced many CHP configurations. Based on these results we formulate the following

\bigskip

\begin{center}
{$\mathbf{Conjecture}$}: 
\end{center}

the set $\left\{ \xi_1, \dots, \xi_k \right\}$ is a {\sl permutation} of the set 
$\left\{ \varphi_1 + \pi/3, \dots, \varphi_k + \pi/3\right\}$

\bigskip

By accepting this conjecture we obtain the following results
\begin{itemize}
\item The set of angles obtained applying the reflection (\ref{eq_refl}) to the set $\{ \xi_1, \dots, \xi_k\}$ corresponds to a particular permutation of the original set, in which the lexicographic order is reversed, e.g. 
\begin{equation}
abc \dots w \rightarrow w \dots cba \nonumber
\end{equation}

\item The $k$ angles $\xi$ may be grouped into $\ell \leq k$ sets, $\{s_1,s_2,\dots, s_\ell\}$, where  each set contains only angles of the same kind; let $n_i$ be the number of elements in each set so that  $n_1 + n_2 + \dots + n_\ell = k$ (we refer to $\left\{ n_i\right\}$ as degeneracies).

\item Let $n_V$ be the number of occupied vertices in the fundamental domain (by occupied vertex we mean that a disk is placed exactly at a vertex);
if $n_V \geq 2$ for a given configuration, a rotation that brings a different occupied vertex in the position $P_1$ produces a degenerate configuration
with a different lexicographic sequence;

\item in special cases the reflection (\ref{eq_refl}) produces one of the configurations obtained with the rotations described at the previous point; to take into account
this behavior we then define the parameter
\begin{equation}
\eta \equiv \left\{ 
\begin{array}{ccc}
1 & , & {\rm reflection \ produces \ a \ configuration } \\
 &   & {\rm already \ obtained \ by \ rotation } \\
2 & , & {\rm otherwise} \\
\end{array}
\right.
\end{equation}

\item  The number of nonequivalent CHP configurations is then simply
\begin{equation}
\# = \max \left( 1,\frac{k!}{\eta \ n_V \ \left(\prod_{i=1}^{\ell} n_i! \right)} \right) \ ,
\label{n_ind}
\end{equation}
where the factor of $\eta$ eliminates the redundant configurations obtained by reflection. Note that the $\max$ is needed only for the case $k=1$, where there is a single configuration.

For the special case of a circle, we have $n_i = 1$, for $1 \leq i \leq k$, $n_V=k$ (all disks are placed at the vertices of a regular polygon of $6k$ sides) and $\eta=2$; 
eq.~(\ref{n_ind}) then reduces to the expression given by Graham and Lubachevsky \cite{Graham97}
\begin{equation} 
\# = \max\left( 1, \frac{(k-1)!}{2} \right) \ .
\end{equation}

\end{itemize}

In addition to the above properties, our numerical results suggest that the angles $\phi$ and $\xi$ obey the relations:
\begin{equation}
\begin{split}
\sum_{j=1}^k \varphi_j &=  k \pi \left( \frac{1}{6}- \frac{1}{\sigma}\right) \\
\sum_{j=1}^k \xi_j &=  k \pi \left( \frac{1}{2}- \frac{1}{\sigma}\right) \\
\end{split} \ .
\nonumber
\end{equation}

\section{Deterministic algorithm}
\label{sec_detalg}

\begin{figure}
\begin{center}
\bigskip\bigskip\bigskip
\includegraphics[width=6cm]{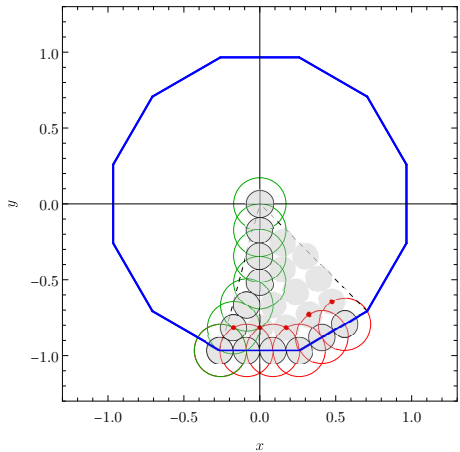}
\bigskip
\caption{Deterministic algorithm for curved hexagonal packing inside a dodecagon. The region enclosed between the dashed lines is the fundamental domain on which we work.}
\label{Fig_Alg_Det}
\end{center}
\end{figure}

Using the results of the previous section we have devised a deterministic algorithm for CHP, that allows one to build in a systematic way {\sl all} the configurations for a given $\sigma$ and
a given $k$. Fig.~\ref{Fig_Alg_Det} illustrates this process for the case of the dodecagon ($\sigma=12$) with $k=6$ ($N=127$).

The algorithm consists of the following steps:
\begin{itemize}
\item[1)] identify the border points of the CHP (these points are {\sl unique} for a given $\sigma$ and $k$) and the corresponding angles, $\varphi_1, \dots, \varphi_k$;
\item[2)] pick a particular permutation of $\left\{ \varphi_1 + \pi/3, \dots, \varphi_k + \pi/3\right\}$, which corresponds to identifying the angles $\{ \xi_1, \dots, \xi_k\}$ 
(we refer to this sequence as to the {\sl DNA} of the configuration because it completely determines the configuration);
\item[3)] as a result of step 2, the first disk of  each internal shell is placed; starting from the outermost shell and moving from left to right,
one can determine the places at which the new disks should be placed by determining the intersections  between the circles centered at the disks that are already placed;
once a new disk is placed, new intersections can be calculated and used to place more disks, until the shell is completed (remember that it is sufficient to work in a region of 
angular width $\pi/3$, because of the symmetry);
\item[4)] every time that a shell is completed move to the next shell and repeat the process; iterate until all the shells are filled;
\end{itemize}

In this way the same configurations that were obtained with the numerical algorithms are now obtained systematically, much more rapidly and with much greater precision. By knowing the "building blocks" $\left\{ \xi_i\right\}$ one is able to calculate {\sl all} the CHP configurations, since they correspond to all possible inequivalent permutations of the building blocks.

\newpage
\section{Results}
\label{sec:res}

In this section we present the results obtained either using the numerical algorithms described in section \ref{sec:algo} or the deterministic algorithm described in section \ref{sec:chp}. We briefly mention also the hexagon and the circle, which have been studied before. In particular curved hexagonal packing was first observed by Graham and Lubachevsky~\cite{Graham97}  for the circle.

\subsection{The hexagon ($\sigma=6$)}
\label{sec:hex}

The hexagon has been recently studied in \cite{Amore22, Amore22hex}.
For $N(k) = 3 k (k+1)+1$ and $k=1,2,\dots$ the best packing configuration  of $N$ congruent disks inside an hexagon corresponds to a perfect hexagonal packing and it is a global maximum of the density (see also \cite{Graham97}). The density (packing fraction)  in this case is given by eq.~(\ref{eq_rho_chp_hex}) and
and $\lim_{k \rightarrow \infty} \rho(k,6) = \rho_{plane} \equiv \frac{\pi}{\sqrt{12}}$. 

\begin{figure}
\begin{center}
\bigskip\bigskip\bigskip
\includegraphics[width=4cm]{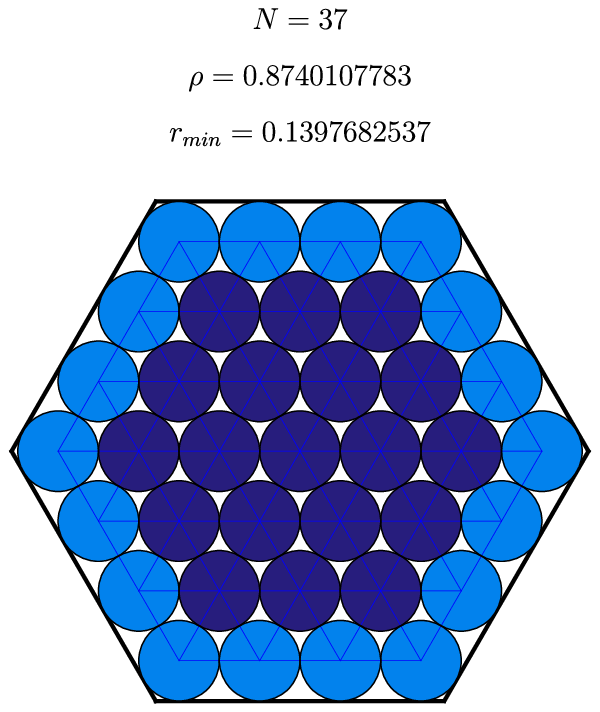}\hspace{1cm} 
\includegraphics[width=4cm]{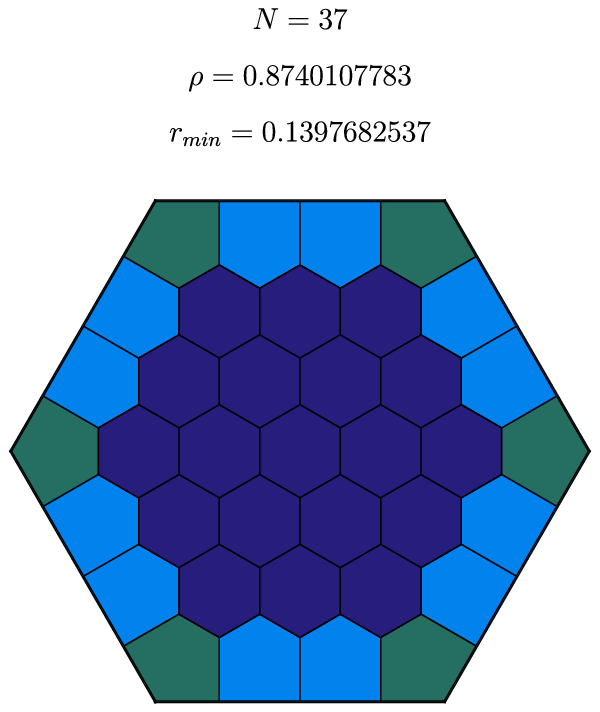}\hspace{1cm} 
\bigskip
\caption{Hexagonal configuration of $37$ congruent disks inside a regular hexagon; the right plot is the Voronoi diagram of the configuration. The colors in these diagrams are assigned depending the number of contacts (left plot) or the sides of the Voronoi cell (right plot).}
\label{Fig_hex}
\end{center}
\end{figure}

For these configurations the internal cells of the Voronoi diagrams are  hexagonal, whereas the border cells are either pentagonal or quadrilateral (see Fig.~\ref{Fig_hex}). 

\subsection{The circle ($\sigma=\infty$)}  

For the circle Graham and Lubachevsky \cite{Graham97}  have found that configurations with curved hexagonal packing exist for the same values of $N(k)$ of the hexagon, with $N \leq 91$ (for $N >91$ these configurations while still existing are not the densest). In eq.(4) of their paper  they provide an explicit expression for the density  and they find that $\lim_{N \rightarrow \infty} \rho^{(circle)}_{\rm CHP}(N) = \frac{\pi^2}{12}$. 

Since $\frac{\pi^2}{12} < \rho_{\rm plane}$ CHP configurations are bound to become non-optimal if $N(k)$ is large enough. The numerical results of \cite{Graham97} identified $N=91$ as the maximal number of disks for which CHP is still optimal.

They also found that for $N(k) = 3 k (k+1)+1$ disks with $k=1,2,\dots$  there are $\max((k-1)!/2,1)$ nonequivalent CHP configurations.

The characterization of the configurations by means of their DNA allow us not only to confirm the results of Graham and Lubachevsky but also explicitly build with extreme ease all the configurations corresponding to a give number of shells. We have included the plots of all configurations up to $7$ shells in the supplemental material \cite{ACZ23a}. 

\subsection{Regular polygons with $\sigma = 6 j$ }

Regular polygons with $6j$ sides ($j=2,3,\dots$) are the natural candidates for possessing configurations with curved hexagonal packing since
the border is invariant under rotations of $\pi/3$. For large enough $k$, however, these configurations are necessarily non optimal since using  eq.(\ref{eq_rho_chp_full})  we obtain
\begin{equation}
\lim_{k \rightarrow \infty} \rho_{\rm CHP} =  \frac{\pi  \sigma}{12}   \tan \left(\frac{\pi }{\sigma }\right) \leq \rho_{\rm plane} \equiv \frac{\pi}{\sqrt{12}} \ ,
\end{equation}
where $\sigma=6$ is the only case where the expression above becomes an equality.

In Tables \ref{Table 1} and \ref{Table 2} we provide a complete classification of the CHP configurations for several regular polygons ($\sigma=12,18,\dots, 60$). To start with, in the eighth and ninth columns we can compare the number of inequivalent configurations obtained with our deterministic algorithm with the explicit formula that we have derived, eq.~(\ref{n_ind}).
Of particular importance is the fourth column, that reports the degeneracies: for all case in the tables we observe an initial sequence with $\left\{n_i\right\} = \left\{ 1,1, \dots,1\right\}$ for $1 \leq k <\sigma/3$ (for $\sigma \geq 30$ however we have not reached sufficiently large systems): this pattern is the same as the circle and corresponds to reaching a denser CHP  since all vertices are occupied.

All the configurations reported in these tables can be found in the supplementary material~\cite{ACZ23a,ACZ23b,ACZ23c,ACZ23d}.

\begin{sidewaystable}
\caption{Independent configurations for the regular polygons with $\sigma=12,18,24,30$ and different number of shells.}
\begin{center}
\begin{tabular}{|c|c|c|c|c|c|c|c|c|}
    \hline
    \multirow{2}{1em}{$\sigma$} & \multirow{2}{1em}{k} & building & \multirow{2}{2em}{$\left\{ n_i\right\}$} & \multirow{2}{1em}{$\eta$} & \multirow{2}{1em}{$n_v$} & \multirow{2}{5em}{Mod(k, $\frac{\sigma}{6}$)} & independent  & \multirow{2}{10em}{$\max \left( 1, \frac{k!}{\eta \ n_v \ \left( \prod_{i} n_i! \right)} \right)$}\\
    & & blocks &  &  & & & configurations &\\
\hline\multirow{8}{1em}{12}
 	 & 1 & 1 & 1 & 1 & 1 & 1 & 1 & 1 \\ 
 	 & 2 & 2 & 1,1 & 1 & 2 & 0 & 1 & 1 \\ 
 	 & 3 & 3 & 1,1,1 & 2 & 1 & 1 & 3 & 3 \\ 
 	 & 4 & 2 & 2,2 & 1 & 2 & 0 & 3 & 3 \\ 
 	 & 5 & 3 & 2,1,2 & 2 & 1 & 1 & 15 & 15 \\ 
 	 & 6 & 2 & 3,3 & 1 & 2 & 0 & 10 & 10 \\ 
 	 & 7 & 3 & 3,1,3 & 2 & 1 & 1 & 70 & 70 \\ 
 	 & 8 & 2 & 4,4 & 1 & 2 & 0 & 35 & 35 \\ 
 	 & 9 & 3 & 4,1,4 & 2 & 1 & 1 & 315 & 315 \\ 
 	 & 10 & 2 & 5,5 & 1 & 2 & 0 & 126 & 126 \\ 
 	 \hline
\multirow{8}{1em}{18}
 	 & 1 & 1 & 1 & 1 & 1 & 1 & 1 & 1 \\ 
 	 & 2 & 2 & 1,1 & 2 & 1 & 2 & 1 & 1 \\ 
 	 & 3 & 3 & 1,1,1 & 2 & 3 & 0 & 1 & 1 \\ 
 	 & 4 & 4 & 1,1,1,1 & 2 & 1 & 1 & 12 & 12 \\ 
 	 & 5 & 5 & 1,1,1,1,1 & 2 & 1 & 2 & 60 & 60 \\ 
 	 & 6 & 3 & 2,2,2 & 2 & 3 & 0 & 15 & 15 \\ 
 	 & 7 & 5 & 2,1,1,1,2 & 2 & 1 & 1 & 630 & 630 \\ 
 	 & 8 & 5 & 2,1,2,1,2 & 2 & 1 & 2 & 2520 & 2520 \\ 
 	 \hline
\multirow{8}{1em}{24}
 	 & 1 & 1 & 1 & 1 & 1 & 1 & 1 & 1 \\ 
 	 & 2 & 2 & 1,1 & 1 & 2 & 2 & 1 & 1 \\ 
 	 & 3 & 3 & 1,1,1 & 2 & 1 & 3 & 3 & 3 \\ 
 	 & 4 & 4 & 1,1,1,1 & 2 & 4 & 0 & 3 & 3 \\ 
 	 & 5 & 5 & 1,1,1,1,1 & 2 & 1 & 1 & 60 & 60 \\ 
 	 & 6 & 6 & 1,1,1,1,1,1 & 2 & 2 & 2 & 180 & 180 \\ 
 	 & 7 & 7 & 1,1,1,1,1,1,1 & 2 & 1 & 3 & 2520 & 2520 \\ 
 	 & 8 & 4 & 2,2,2,2 & 2 & 4 & 0 & 315 & 315 \\ 
 	 \hline
\multirow{8}{1em}{30}
 	 & 1 & 1 & 1 & 1 & 1 & 1 & 1 & 1 \\ 
 	 & 2 & 2 & 1,1 & 2 & 1 & 2 & 1 & 1 \\ 
 	 & 3 & 3 & 1,1,1 & 2 & 1 & 3 & 3 & 3 \\ 
 	 & 4 & 4 & 1,1,1,1 & 2 & 1 & 4 & 12 & 12 \\ 
 	 & 5 & 5 & 1,1,1,1,1 & 2 & 5 & 0 & 12 & 12 \\ 
 	 & 6 & 6 & 1,1,1,1,1,1 & 2 & 1 & 1 & 360 & 360 \\ 
 	 & 7 & 7 & 1,1,1,1,1,1,1 & 2 & 1 & 2 & 2520 & 2520 \\ 
 	 & 8 & 8 & 1,1,1,1,1,1,1,1 & 2 & 1 & 3 & 20160 & 20160 \\ 
 	 \hline
   \end{tabular}
\end{center}
\label{Table 1}
\end{sidewaystable}

\begin{sidewaystable}
\caption{Independent configurations for the regular polygons with $\sigma=36,42,48,54,60$ and different number of shells.}
\begin{center}
\begin{tabular}{|c|c|c|c|c|c|c|c|c|}
    \hline
    \multirow{2}{1em}{$\sigma$} & \multirow{2}{1em}{k} & building & \multirow{2}{2em}{$\left\{ n_i\right\}$} & \multirow{2}{1em}{$\eta$} & \multirow{2}{1em}{$n_v$} & \multirow{2}{5em}{Mod(k, $\frac{\sigma}{6}$)} & independent  & \multirow{2}{10em}{$\max \left( 1, \frac{k!}{\eta \ n_v \ \left( \prod_{i} n_i! \right)} \right)$}\\
    & & blocks &  &  & & & configurations &\\
    \hline
\multirow{8}{1em}{36}
 	 & 1 & 1 & 1 & 1 & 1 & 1 & 1 & 1 \\ 
 	 & 2 & 2 & 1,1 & 1 & 2 & 2 & 1 & 1 \\ 
 	 & 3 & 3 & 1,1,1 & 2 & 3 & 3 & 1 & 1 \\ 
 	 & 4 & 4 & 1,1,1,1 & 2 & 2 & 4 & 6 & 6 \\ 
 	 & 5 & 5 & 1,1,1,1,1 & 2 & 1 & 5 & 60 & 60 \\ 
 	 & 6 & 6 & 1,1,1,1,1,1 & 2 & 6 & 0 & 60 & 60 \\ 
 	 & 7 & 7 & 1,1,1,1,1,1,1 & 2 & 1 & 1 & 2520 & 2520 \\ 
 	 & 8 & 8 & 1,1,1,1,1,1,1,1 & 2 & 2 & 2 & 10080 & 10080 \\ 
 	 \hline
\multirow{8}{1em}{42}
 	 & 1 & 1 & 1 & 1 & 1 & 1 & 1 & 1 \\ 
 	 & 2 & 2 & 1,1 & 2 & 1 & 2 & 1 & 1 \\ 
 	 & 3 & 3 & 1,1,1 & 2 & 1 & 3 & 3 & 3 \\ 
 	 & 4 & 4 & 1,1,1,1 & 2 & 1 & 4 & 12 & 12 \\ 
 	 & 5 & 5 & 1,1,1,1,1 & 2 & 1 & 5 & 60 & 60 \\ 
 	 & 6 & 6 & 1,1,1,1,1,1 & 2 & 1 & 6 & 360 & 360 \\ 
 	 & 7 & 7 & 1,1,1,1,1,1,1 & 2 & 7 & 0 & 360 & 360 \\ 
 	 & 8 & 8 & 1,1,1,1,1,1,1,1 & 2 & 1 & 1 & 20160 & 20160 \\ 
 	 \hline
\multirow{8}{1em}{48}
 	 & 1 & 1 & 1 & 1 & 1 & 1 & 1 & 1 \\ 
 	 & 2 & 2 & 1,1 & 1 & 2 & 2 & 1 & 1 \\ 
 	 & 3 & 3 & 1,1,1 & 2 & 1 & 3 & 3 & 3 \\ 
 	 & 4 & 4 & 1,1,1,1 & 2 & 4 & 4 & 3 & 3 \\ 
 	 & 5 & 5 & 1,1,1,1,1 & 2 & 1 & 5 & 60 & 60 \\ 
 	 & 6 & 6 & 1,1,1,1,1,1 & 2 & 2 & 6 & 180 & 180 \\ 
 	 & 7 & 7 & 1,1,1,1,1,1,1 & 2 & 1 & 7 & 2520 & 2520 \\ 
 	 & 8 & 8 & 1,1,1,1,1,1,1,1 & 2 & 8 & 0 & 2520 & 2520 \\ 
 	 \hline
\multirow{8}{1em}{54}
 	 & 1 & 1 & 1 & 1 & 1 & 1 & 1 & 1 \\ 
 	 & 2 & 2 & 1,1 & 2 & 1 & 2 & 1 & 1 \\ 
 	 & 3 & 3 & 1,1,1 & 2 & 3 & 3 & 1 & 1 \\ 
 	 & 4 & 4 & 1,1,1,1 & 2 & 1 & 4 & 12 & 12 \\ 
 	 & 5 & 5 & 1,1,1,1,1 & 2 & 1 & 5 & 60 & 60 \\ 
 	 & 6 & 6 & 1,1,1,1,1,1 & 2 & 3 & 6 & 120 & 120 \\ 
 	 & 7 & 7 & 1,1,1,1,1,1,1 & 2 & 1 & 7 & 2520 & 2520 \\ 
 	 & 8 & 8 & 1,1,1,1,1,1,1,1 & 2 & 1 & 8 & 20160 & 20160 \\ 
 	 \hline
\multirow{8}{1em}{60}
 	 & 1 & 1 & 1 & 1 & 1 & 1 & 1 & 1 \\ 
 	 & 2 & 2 & 1,1 & 1 & 2 & 2 & 1 & 1 \\ 
 	 & 3 & 3 & 1,1,1 & 2 & 1 & 3 & 3 & 3 \\ 
 	 & 4 & 4 & 1,1,1,1 & 2 & 2 & 4 & 6 & 6 \\ 
 	 & 5 & 5 & 1,1,1,1,1 & 2 & 5 & 5 & 12 & 12 \\ 
 	 & 6 & 6 & 1,1,1,1,1,1 & 2 & 2 & 6 & 180 & 180 \\ 
 	 & 7 & 7 & 1,1,1,1,1,1,1 & 2 & 1 & 7 & 2520 & 2520 \\ 
 	 & 8 & 8 & 1,1,1,1,1,1,1,1 & 2 & 2 & 8 & 10080 & 10080 \\ 
 	 \hline
   \end{tabular}
   \end{center}
\label{Table 2}
   \end{sidewaystable}

\begin{table}
\caption{Emergence of non CHP configurations}
\begin{center}
\begin{tabular}{|c|c|c|c|}
\hline
$\sigma$ & $k$ & $\rho$ & $\rho-\rho_{CHP}$ \\
\hline
12 & 7 & 0.838209 & 0.0109501 \\
\hline
18 & 7 & 0.826627 & 0.00687681 \\
\hline
24 & 6 & 0.818231 & 0.000477362 \\
\hline
30 & 6 & 0.817599 & 0.00137266 \\
\hline
36 & 7 & 0.825064 & 0.00798317 \\
\hline
42 & 6 & 0.817919 & 0.0017007 \\
\hline
48 & 6 & 0.816806 & 0.000371802 \\
\hline
54 & 6 & 0.816963 & 0.000313375 \\
\hline
60 & 6 & 0.816443 & 0.0000973674 \\
\hline
66 & 6 & 0.816631 & 0.000417085 \\
\hline
72 & 7 & 0.825047 & 0.00794864 \\
\hline
78 & 6 & 0.816701 & 0.000480921 \\
\hline
84 & 6 & 0.816578 & 0.000294333 \\
\hline
90 & 6 & 0.816688 & 0.000323156 \\
\hline
96 & 6 & 0.816786 & 0.000517179 \\
\hline
\end{tabular}
\end{center}
\label{Table3}
\end{table}

\begin{figure}
\begin{center}
\bigskip\bigskip\bigskip
\includegraphics[width=8cm]{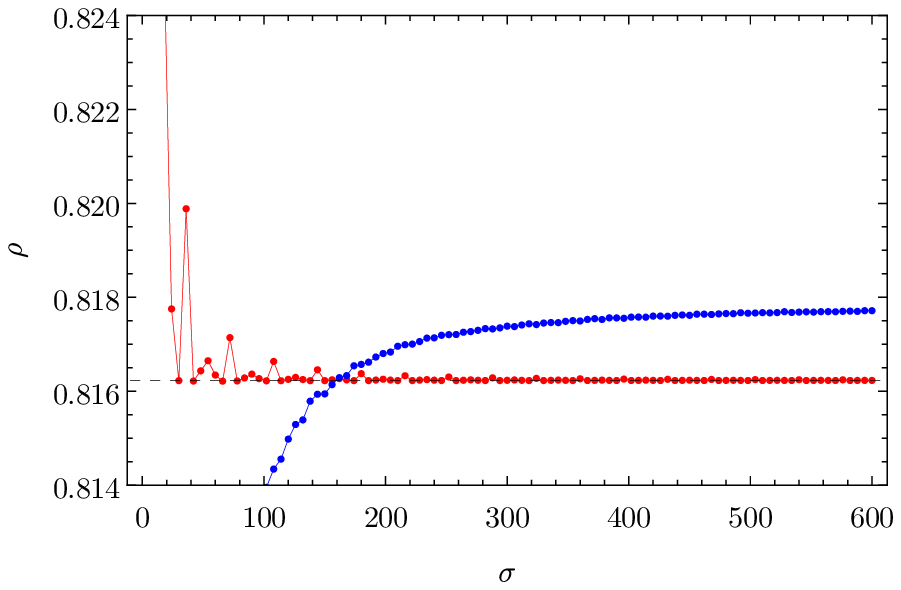}\hspace{1cm} 
\caption{Density of the CHP configurations (red curve) and of the non--CHP configuration (blue curve) for different $\sigma$.}
\label{Fig_break_CHP}
\end{center}
\end{figure}

In Table \ref{Table3}  we report the first cases where non-CHP configurations are found for polygons with $\sigma= 12, 18, \dots, 96$. The cases of the dodecagon, octadecagon, triacontakaihexagon and heptacontakaidigon ($\sigma=12,18,36,72$) stand out, because are the unique where we could not find a configuration denser than CHP for $k=6$ shells ($N=127$); of course, this does not amount to a formal proof that for these polygons the densest configuration of $127$ disks is CHP, although our numerical findings support this conclusion. The remaining polygons in the table, on the other hand, share the same fate of the circle, for which Graham and Lubachevsky showed that $N=91$ is the last case where CHP are densest~\cite{Graham97}.

Intuitively we expect that for large enough  $\sigma$ the CHP will not be optimal, since the densest configuration of $127$ disks inside the circle ($\sigma=\infty$) 
is not CHP~\cite{Graham97,SpechtRepo}.  In Fig.~\ref{Fig_break_CHP} we plot the density of the packing configuration of $N=127$ disks inside a regular polygon of $\sigma= 6j$ sides ($j=1,2,\dots$): the red curve is the density of the CHP configurations, given by formula (\ref{eq_rho_chp}), whereas the blue curve is the density of the configurations obtained taking the densest configuration of $127$ disks inside a circle from \cite{SpechtRepo}, adapting it to the regular polygons~\footnote{The points of the original configuration that are already inside the polygon are left untouched, whereas the points that fall in the circle but outside the polygon are projected on the border of the polygon.}. The black dashed curve is the density of the CHP configuration of the circle given in eq.~(\ref{eq_rho_chp_circle}) for $k=6$.
Of course we expect that the configuration of \cite{SpechtRepo} is a good ansatz only for $\sigma \gg 1$: the fact that the blue curve is below the red curve for 
$\sigma < 100$ does not imply that the CHP there is denser than any other configuration, but rather that the ansatz based on the circle is not good enough for a regular polygon with not so many sides. For $\sigma < 100$, however, we have the numerical results of Table \ref{Table3} which show that most of configurations for $k=6$ are 
indeed non-CHP. From the plot we see that for $\sigma \geq 162$ the modified ansatz is always denser than the corresponding CHP.

\begin{table}
\caption{Number of border disks for the first non--CHP configurations occurring in regular polygons with $6j$ sides.}
\begin{center}
\begin{tabular}{|c|c|c|c||c|c|c|c||c|c|c|c|}
\hline 
$\sigma$ & $k$ & $N_b$ & $6k$ & $\sigma$ & $k$ & $N_b$ & $6k$ & $\sigma$ & $k$ & $N_b$ & $6k$ \\
\hline 
12       & 7  & 35  & 42  & 18  & 7 & 33 & 42 & 24 & 6 & 29 & 36\\
\hline 
30       & 6  & 32  & 36  & 36  & 7 & 38 & 42 & 42 & 6 & 31 & 36\\
\hline 
48       & 6  & 31  & 36  & 54  & 6 & 29 & 36 & 60 & 6 & 32 & 36\\
\hline 
66       & 6  & 30  & 36  & 72  & 7 & 37 & 42 & 78 & 6 & 31 & 36\\
\hline 
84       & 6  & 30  & 36  & 90  & 6 & 31 & 36 & 96 & 6 & 31 & 36 \\
\hline
\end{tabular}
\end{center}
\label{Table4}
\end{table}

In Table \ref{Table4} we report the number of border disks for the first non--CHP configurations that are found to be denser than the corresponding CHP configurations, for the same regular polygons of Table \ref{Table3}. The striking aspect in this table is the fact that $N_b$ is always much smaller than $6k$, the number of border disks in a CHP. For larger $k$ the factor $6k-N_b$ is expected to grow.

The diagrams for the non--CHP configurations reported in Tables~\ref{Table3} and \ref{Table4} are found in \cite{ACZ23e} (\cite{ACZ23f} contains the complete set of data corresponding to the configurations obtained).

\subsection{Exceptional CHP}

By exceptional CHP, we mean packing configurations in the regular polygon which are invariant under $\pi/3$ rotations, but either $\sigma \neq 6j$ (so that the border  {\sl is not invariant under the rotation}) or $N \neq N(k)$, so that the shell structure of the original hexagonal packing is not retained.

\begin{figure}
\begin{center}
\bigskip\bigskip\bigskip
\includegraphics[width=4cm]{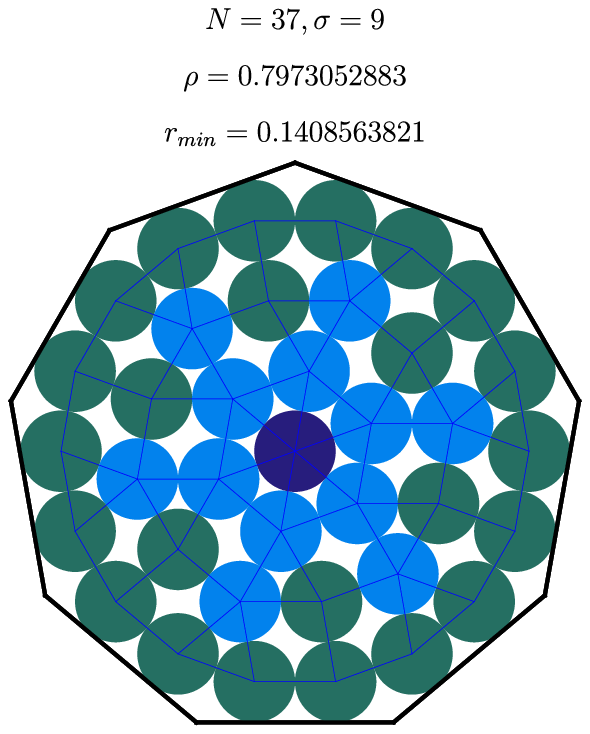}\hspace{1cm} 
\includegraphics[width=4cm]{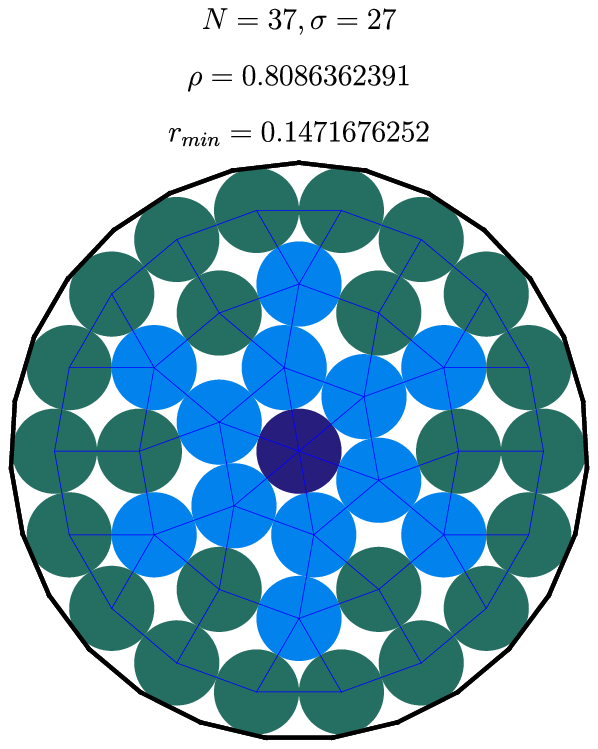}\hspace{1cm} 
\caption{Exceptional hexagonal packing of 37 disks inside a nonagon (left) and in a icosiheptagon (right)}
\label{Fig_EHP}
\end{center}
\end{figure}

In the first class we have found two examples, corresponding to the packing of $37$ disks inside a nonagon or an icosiheptagon (see Fig.~\ref{Fig_EHP}).
It is remarkable that these configurations are  invariant under rotations of $\pi/3$ even though the border is not: as a matter of fact
we have performed several numerical experiments over different polygons and with different (hexagonal) numbers of disks, but we have failed to find additional cases (while this does not mean that there are no other cases, it signals however that these configurations are extremely rare).

\begin{figure}
\begin{center}
\bigskip\bigskip\bigskip
\includegraphics[width=4cm]{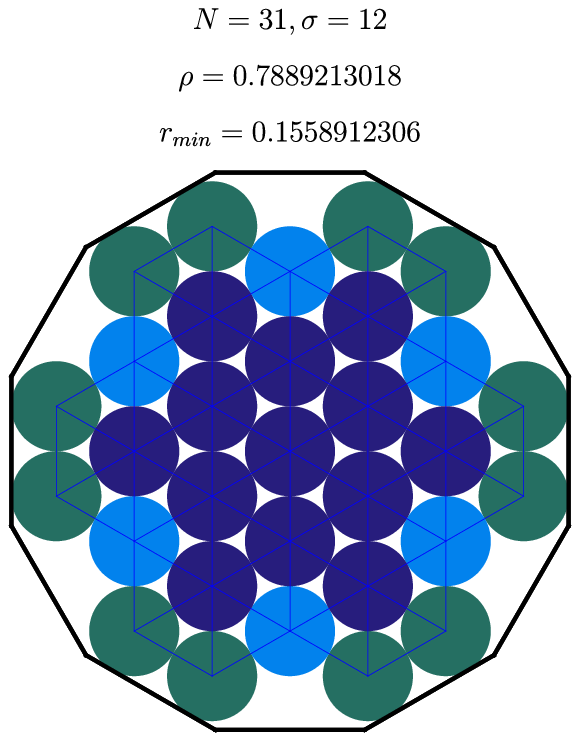}\hspace{1cm} 
\includegraphics[width=4cm]{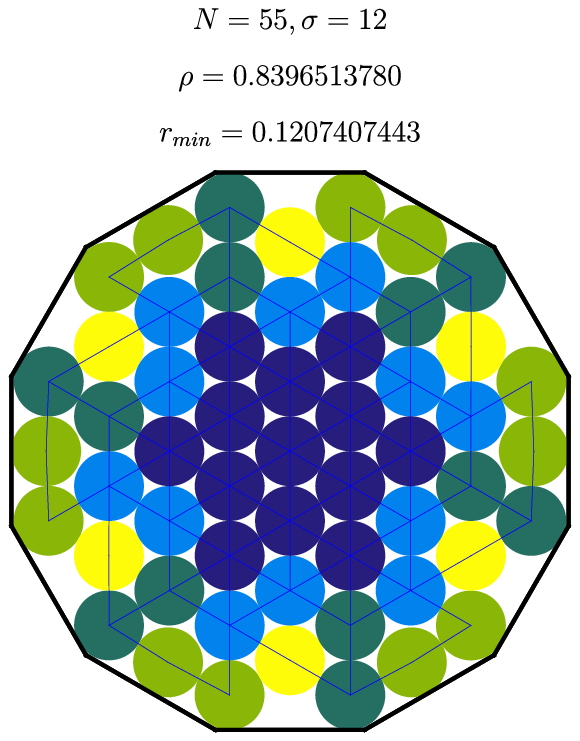}\hspace{1cm} 
\caption{Exceptional hexagonal packing of 31 (left) and 55 (right)  disks inside a dodecagon}
\label{Fig_EHPb}
\end{center}
\end{figure}

In the second class we have found examples in several regular polygons with $6j$ sides, particularly for $N=31$ and $N=55$, which are displayed in Fig.~\ref{Fig_EHPb}  for the dodecagon. In this case the configurations consist of a tightly packed core with perfect hexagonal packing, surrounded by a less dense region where peripheral disks are distributed symmetrically along the border. For the configurations in the figure, the cores correspond to the hexagonal numbers $19$ and $37$, with $12$ and $18$ border disks: on the other hand, the configurations of this kind with a core of $7$ disks and $6$ border disks, or a core with $61$ disks and $24$ border disks {\sl are not} global maxima of the density. This configurations appear to be related to the "wedge hexagonal packing" configurations first identified by Hopkins in \cite{Hopkins12}.

\newpage

\section{Breaking the ice}
\label{sec:break}

In this section we want to  explore a further application of our findings: although configurations with curved hexagonal packing are optimal only up to some critical number of disks $N_{\rm max}$ (for the circle $N_{\rm max}=91$), one can use them as starting point to generate much denser irregular configurations, by applying algorithm 2 of \cite{Amore22} once or repeated times. If $N$ is large enough the system will arrange in a disordered configuration, where the structure of the original crystal is completely washed out.

\begin{figure}
\begin{center}
\bigskip\bigskip\bigskip
\includegraphics[width=8cm]{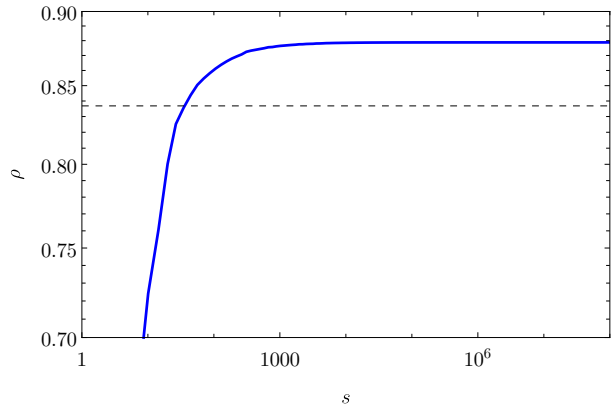}\hspace{1cm} 
\caption{Evolution of the density in a single run of algorithm 2 for $N=1657$ disks inside a dodecagon. The dashed line is the density of the CHP configurations, $\rho_{\rm CHP} \approx 0.8368374943$.}
\label{Fig_ALG2_ICE}
\end{center}
\end{figure}

In fig.~\ref{Fig_ALG2_ICE} we show the evolution of the density of a configuration with $1657$ disks inside a dodecagon, in a single run of algorithm 2: at the beginning of the algorithm a CHP configuration (in our case the one of smallest lexicographic order) is constructed and the positions of the disks are randomly altered by a small amount;
this process typically produces an initial ansatz with a smaller density than the original CHP configuration (the dashed line in the plot). As the algorithm proceeds, the value of $s$ (the exponent in the potential which controls its range of action) is gradually increased, up to very large values, $s=10^8$ in this case. During this process, the configuration starts reaccomodating and gradually increases its density, reaching a much denser configuration (already for $s \approx 40$ the system improves the CHP density). There is no guarantee of course that the new configuration, while much better than the original, will be a global maximum of the packing fraction (most likely the reverse is true!) but the algorithm 2 can be applied iteratively, each time on the best configuration found, thus leading to further improvements of the density.
The focus of our present analysis however is not finding the global maximum, which for such large values of $N$ would be extremely challenging and time consuming, but rather show that it is possible to generate very dense configuration for $N \gg 1$ in a simple and effective way.

\begin{figure}
\begin{center}
\bigskip\bigskip\bigskip
\includegraphics[width=4cm]{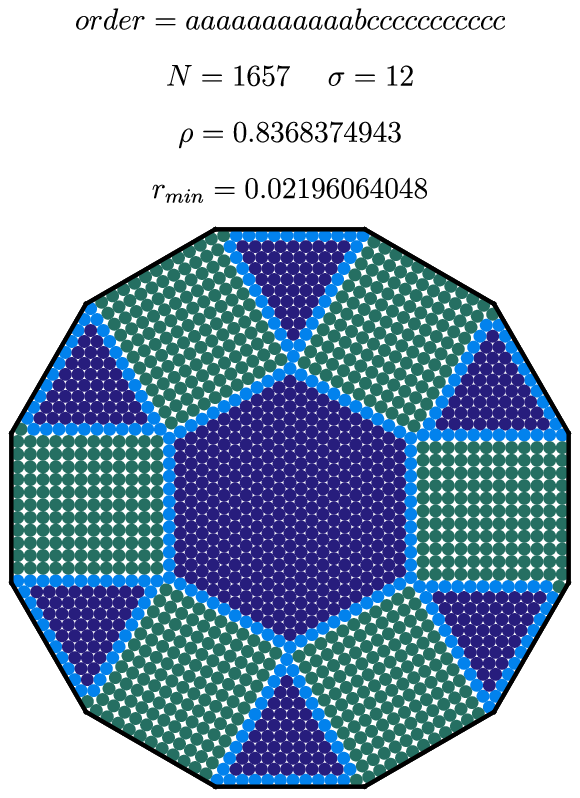}\hspace{1cm} 
\includegraphics[width=4cm]{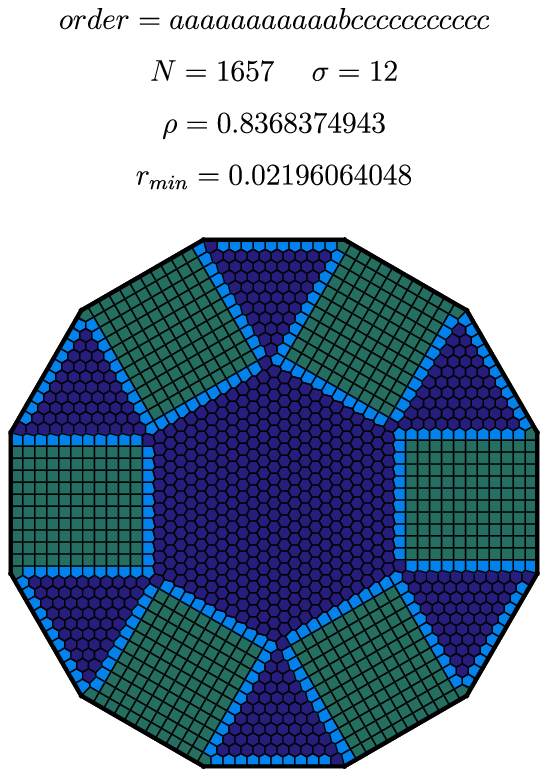} 
\caption{CHP configuration of lowest lexicographic order for  $k=23$ ($N=1657$) inside a dodecagon. The right plot is the Voronoi diagram.}
\label{Fig_CHP_12_23}
\end{center}
\end{figure}

\begin{figure}
\begin{center}
\bigskip\bigskip\bigskip
\includegraphics[width=4cm]{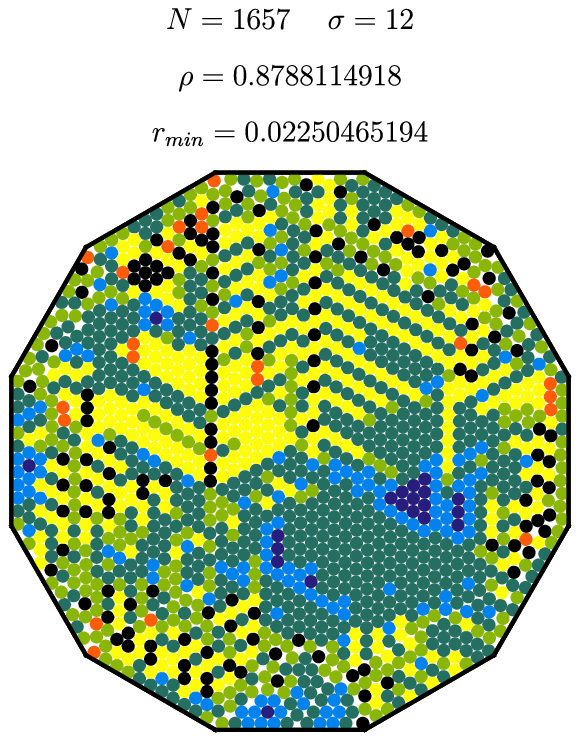}\hspace{1cm} 
\includegraphics[width=4cm]{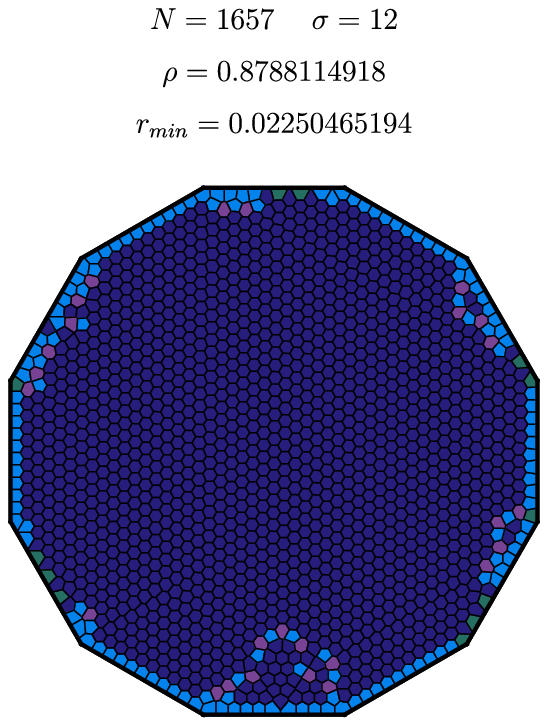} 
\caption{Configuration of $N=1657$ disk inside the dodecagon found with a single run of algorithm 2 to the configuration of Fig.~\ref{Fig_CHP_12_23}.}
\label{Fig_ALG2_1657}
\end{center}
\end{figure}

In Fig.~\ref{Fig_CHP_12_23}  we display the CHP configuration of lowest lexicographic order inside a dodecagon with $N=1657$ disks: the right plot is the Voronoi diagram.
The much denser configuration in Fig.~\ref{Fig_ALG2_1657}  has been obtained  with a single run of algorithm 2, taking as initial ansatz the previous CHP configuration, with a small random perturbation. Although in both cases the topological charge of the Voronoi diagrams must add up to $6$, due to Euler's theorem, in the second case we notice two notable facts: the disappearance of higher order vertices (which carry a negative topological charge) and the almost complete expulsion of the topological charge to the border of the domain. See \cite{Amore22} for a discussion of Euler's theorem and topological charge.

\newpage

\section{Conclusions}
\label{sec:concl}

The main findings of this paper can be resumed as follows:
\begin{itemize}
\item the densest packing configuration of  $N(k) = 3 k (k+1)+1$ congruent disks inside a regular polygon with $\sigma = 6j$ sides corresponds for $k \leq k_{\rm max}$ to a curved hexagonal packing;

\item we have been able to {\sl fully} characterize these configurations and set up a deterministic algorithm that can ensemble {\sl all} the configurations for a given $\sigma$ and $k$;

\item The CHP configurations have a density

\begin{equation}
\begin{split}
\rho(k,\sigma) &= 8 \pi  (3 k (k+1)+1) \frac{ \csc \left(\frac{2 \pi }{\sigma }\right) \left(\sin \left(\frac{\pi }{\sigma }\right) +\cos \left(\frac{\pi}{\sigma }+\frac{\pi }{6}\right)\right)^2}{\sigma  \left(4 \sum_{j=1}^k \cos (\phi (j))+\tan \left(\frac{\pi}{\sigma}\right)+\sqrt{3}\right)^2} 
\end{split}
\nonumber
\end{equation}

\item these configurations appear to become highly degenerate (in the density) as $k$ gets larger; 
the number of nonequivalent CHP configurations is
\begin{equation}
\# = \max \left( 1, \frac{k!}{\eta \ n_v \ \left( \prod_{i} n_i! \right)} \right) \nonumber
\end{equation}

\item we have found curved hexagonal packing also to be present for selected regular polygons of $3 (2j+1)$ sides, for specific number of disks, even though the border of the polygon is invariant under rotations of $2\pi/3$, instead of $\pi/3$;

\item there are special configurations, at selected values of $N$, in which a perfect hexagonal packing is present in the inner region, but is disrupted at the border; 

\item our analysis confirms the main results of Graham and Lubachevsky for curved hexagonal packing in the circle, (we reproduce both the density and the number of configurations reported in ref.~\cite{Graham97}), but additionally we provide a simple unified picture of all curved hexagonal packing in regular polygons (including the circle) and a powerful deterministic algorithm that allows one to explicitly calculate all the configurations corresponding to a given number of shells;

\item we show that it is possible to generate very dense disordered systems with very large number of constituents, by using a CHP configuration as initial ansatz and applying the algorithm 2 of \cite{Amore22, Amore22b} to it.

\end{itemize}

\section*{Acknowledgements}
The research of P.A. was supported by Sistema Nacional de Investigadores (M\'exico). 
The plots in this paper have been plotted using Mathematica~\cite{wolfram} and {\rm MaTeX} \cite{szhorvat}. 
Numerical calculations have been carried out using python ~\cite{python} and numba ~\cite{numba}.

\end{document}